\documentclass[11pt]{amsart}   	
\pdfoutput=1

\usepackage[margin=1in]{geometry}
					
\usepackage{amsmath,amssymb,amsthm}
\usepackage{stmaryrd}
\usepackage{MnSymbol}
\usepackage{hyperref}
\usepackage{mathtools}
\usepackage[capitalise]{cleveref}
\usepackage{tikz-cd}
\usepackage{tikz}
\usepackage{enumitem}
\usepackage{makecell}
\setcellgapes{4pt}
\usepackage{pdfpages} 
\usepackage{enumitem}
\usepackage{subcaption}
\usepackage{multicol}
\usepackage{pgfplots, pgfplotstable}

\usepackage{enumitem}
\usepackage{tabstackengine}

\usepackage{comment}
\usepackage{pdflscape}

\newcommand{\add}{\mathsf{add}}

\newcommand{\mods}{\mathsf{mod}}

\newcommand{\proj}{\mathsf{proj}}

\newcommand{\cone}{\mathsf{Cone}}
\newcommand{\Hasse}{\mathsf{Hasse}}

\newcommand{\N}{\mathbb{N}}

\newcommand{\K}{\mathrm{K}}

\newcommand{\mfrak}{\mathfrak{a}}

\DeclareMathOperator{\Hom}{\mathrm{Hom}}

\DeclareMathOperator{\End}{\mathrm{End}}

\DeclareMathOperator{\rad}{\mathrm{rad}}

\DeclareMathOperator{\thick}{\mathrm{thick}}
\DeclareMathOperator{\silt}{\mathrm{silt}}

\DeclareMathOperator{\twosilt}{\mathrm{2-silt}}

\newtheorem{theorem}{Theorem}[section]
\newtheorem{corollary}[theorem]{Corollary}
\newtheorem{lemma}[theorem]{Lemma}

\newtheorem{proposition}[theorem]{Proposition}

\theoremstyle{definition}
\newtheorem{definition}[theorem]{Definition}

\newtheorem{example}[theorem]{Example}

\usepackage[style=alphabetic, maxnames=1000]{biblatex}

\title{On silting theory for complete path algebras}

\author[M. Kaipel]{Maximilian Kaipel}
\address[MK]{Fakultät für Mathematik, Universität Bielefeld, 33501 Bielefeld, Germany}
\email{mkaipel@math.uni-bielefeld.de}
\thanks{MK acknowledges support from the Deutsche Forschungsgemeinschaft (DFG, German Research Foundation) as part of SFB-TRR 191 (Project ID 281071066) and as part of SFB-TRR 358/1 2023 (Project ID 491392403).}

\author[D. Labardini-Fragoso]{Daniel Labardini-Fragoso}
\address[DLF]{Dipartimento di Matematica 'Tullio Levi-Civita', Università degli Studi di Padova, Via Trieste, 63 - 35121 Padova, Italy}
\email{daniel.labardinifragoso@unipd.it}
\thanks{All authors received financial support of Fondazione Cariparo through a “Starting Package per attrarre ricercatrici e ricercatori eccellenti dall’estero” (CUP:
C93C22008360007) granted to DLF at Università degli Studi di Padova}

\author[H. U. Terland]{Håvard Utne Terland}
\address[HUT]{Department of Mathematical Sciences, Norwegian University of Science and Technology (NTNU), 7491 Trondheim, Norway}
\email{havard.u.terland@ntnu.no}
\thanks{}

\begin{document}

\begin{abstract}
	Let $Q$ be a finite quiver and let $\Lambda$ be the complete path algebra of $Q$ with respect to the arrow ideal. We study the mutation of silting objects in the bounded homotopy category of projective $\Lambda$-modules and characterise in which cases the mutation theories of silting complexes and 2-term silting complexes of $\Lambda$ are determined by those of a finite dimensional quotient $\Lambda/I$. This characterisation is closely related to a characterisation of 2-term silting finite quotient algebras. 
\end{abstract}

\maketitle

\section{Introduction}

Silting objects and silting subcategories of triangulated categories originate in the study of $t$-structures on triangulated categories \cite{KV88}. In contrast to tilting objects \cite{Rickard1989}, silting objects do not induce derived equivalences, yet they still encode important homological information. For example, silting objects in the bounded homotopy category of a finite dimensional algebra classify bounded $t$-structures with length heart in the bounded derived category of modules \cite{KY14}. Remarkably, silting subcategories often admit a well-behaved mutation theory \cite{Aihara_2012}, which establishes a procedure for obtaining new silting subcategories from an initial one. 

In representation theory, the triangulated categories of interest are usually Krull--Schmidt categories, which consequently means that it is equivalent to study silting objects instead of silting subcategories. In particular, this holds for the bounded homotopy categories of finitely generated projective modules of finite dimensional algebras and of complete path algebras of quivers \cite{KellerYang2011}. Over finite dimensional algebras, mutation of silting objects is always possible, while over complete path algebras, mutation is not always be possible. The aim of this paper is to describe this apparent discrepancy.

On the one hand, path algebras of quivers are classical objects in representation theory. Gabriel's theorem \cite{Gabriel1972} classifies representation finite path algebras via quivers whose underlying graph is of Dynkin type ADE. Moreover, every algebra over an algebraically closed field is Morita equivalent to a quotient of a path algebra. On the other hand, the significance of complete path algebras has been highlighted by the introduction of quivers with potential \cite{DWZ2008, DWZ2010}. The main result of this paper is a comparison between the silting theory of complete path algebras and the silting theory of their finite dimensional quotient algebras, focusing especially on their behaviour under mutation. This is inspired by the following similar results:
\begin{itemize}
	\item Complete group algebras of pro-finite groups are, just like complete path algebras, examples of pseudo-compact (topological) algebras. This is a consequence of the fact that, they are constructed as inverse limits of group algebras of finite quotient groups. The global dimensions of complete group algebras of pro-finite groups have been described using the global dimensions of the complete group algebras of their Sylow subgroups \cite{Brumer1966}.
	\item $R$-algebras which are finitely generated as a module over a commutative complete local Noetherian ring $R$ are generally infinite dimensional over the residue field of $R$. Surprisingly, however, their silting theory \cite{Gnedin2022}, and in particular their 2-term silting theory \cite{Kimura2024}, is governed by finite dimensional quotient algebras. Particularly important examples of these are finite dimensional $K$-algebras, where $K$ is a field; $R$-orders, such as the group algebra $RG$, where $G$ is a finite group (for which a special instance of Kimura's theorem was shown in \cite{EJR2018}); and complete string algebras in the sense of \cite{Ricke}. These are examples of string algebras over local rings as in \cite{BT2025} for which many exotic examples exist \cite{BT2024}.
\end{itemize}

Silting objects admit a partial ordering $\leq$ determined by the vanishing of $\Hom$ in positive direction \cite{Aihara_2012}. However, it is extremely difficult to describe this partially ordered set without using irreducible mutation as a tool for generating it. For this reason, we focus on the subposet $\leq_{\textnormal{mut}}$ induced by irreducible mutation (\cref{lem:submutationposet}). Over finite dimensional algebras, and so-called Noetherian $R$-algebras as above, these posets coincide but for complete path algebras the posets are generally different. This subtlety is addressed by the following characterisation.

\begin{theorem}\label{thm:main1}
    Let $Q$ be a finite connected quiver and let $\Lambda = \widehat{KQ}$ be its complete path algebra with respect to the arrow ideal. Then the following are equivalent:
    \begin{enumerate}
    \item\label{thm:main1eq1} There exists an ideal $I \subseteq \rad \Lambda$ of $\Lambda$ such that $\Lambda/I$ is finite dimensional and such that there exists an isomorphism of posets
    \[ f: (\twosilt \Lambda, \leq_{\textnormal{mut}}) \simeq (\twosilt \Lambda/I, \leq). \]
    \item\label{thm:main1eq2} There exists an ideal $I \subseteq \rad \Lambda$ of $\Lambda$ such that $\Lambda/I$ is finite dimensional and such that there exists an isomorphism of posets
    \[ f': (\silt \Lambda, \leq_{\textnormal{mut}}) \simeq (\silt \Lambda/I, \leq). \]
    \item\label{thm:main1eq3} $Q$ satisfies one of the following:
    \begin{enumerate}
        \item $Q$ has only one vertex;
        \item $Q$ is acyclic;
        \item $Q$ is an oriented cycle.
    \end{enumerate}
    \item\label{thm:main1eq4} $\Lambda$ is one of the following:
    \begin{enumerate}
    	\item a local $K$-algebra;
	\item a finite dimensional $K$-algebra;
	\item a hereditary $K \llbracket x \rrbracket$-order.
    \end{enumerate}   
    \end{enumerate}
\end{theorem}

We highlight the following ring theoretic corollary.

\begin{corollary}\label{cor:main1}
	Let $Q$ be a finite connected quiver with more than one vertex and let $\Lambda = \widehat{KQ}$ be its complete path algebra with respect to the arrow ideal. Then the following are equivalent:
	\begin{enumerate}
	\item\label{cor:main1eq1} $\Lambda$ is a module-finite $R$-algebra for some commutative complete local Noetherian ring $R$;
	\item\label{cor:main1eq2}  $Q$ is an oriented cycle or acyclic. 
	\end{enumerate}
\end{corollary}

Aside from the pathological silting theory for local algebras \cite[Thm. 2.26]{Aihara_2012}, \cref{thm:main1} shows that it is possible to describe the silting theory of a complete path algebra $\Lambda$ as that of a finite dimensional quotient only if $\Lambda$ is a Noetherian $R$-algebra in the sense of \cite{Ausl1977}. It is shown in \cite{Gnedin2022, Kimura2024} that it is possible to describe the silting theory of a Noetherian $R$-algebra, with an appropriate choice of $R$, in terms of a finite dimensional quotient algebra and we use their results to establish the implications \eqref{thm:main1eq4}$\Rightarrow$\eqref{thm:main1eq1} and \eqref{thm:main1eq4}$\Rightarrow$\eqref{thm:main1eq2} in \cref{thm:main1}.

The restriction of silting theory to complexes of projective modules concentrated in cohomological degrees $-1$ and $0$ has become a central research theme, known as 2-term silting theory, due to its connection to $\tau$-tilting theory \cite{tau, DerksenFei2015} and cluster-tilting theory \cite{BMRRT}.  This theory is especially powerful when only finitely many 2-term silting objects exist. For Noetherian $R$-algebras, 2-term silting objects provide significant insights into the torsion classes in the category of modules  in this case \cite[Thm. 4.10]{IyamaKimura2024}. Motivated by this, we classify for which quivers all finite dimensional quotients have finitely many 2-term silting complexes.

\begin{proposition}\label{thm:main2}
    Let $Q$ be a finite connected quiver and let $\Lambda = \widehat{KQ}$ be the complete path algebra with respect to the arrow ideal $\mathfrak{a}$. Then the following properties are equivalent:
    \begin{enumerate}
    	\item Every finite dimensional quotient of $\Lambda$ is 2-term silting finite;
	\item The quotients of the form $\Lambda/\mfrak^i$ are 2-term silting finite for all $i \in \N$;
	\item $Q$ satisfies one of the following:
	\begin{enumerate}
       		\item$Q$ has only one vertex;
        		\item $Q$ is of Dynkin type ADE;
        		\item $Q$ is an oriented cycle.
        \end{enumerate}
    \end{enumerate}
\end{proposition}

We note that there is a large intersection between the quivers appearing in \cref{thm:main1} and those appearing in \cref{thm:main2}. Using some of the key steps in the proof of \cref{thm:main1}, such as \cref{lem:mutprobcyclein} and \cref{lem:mutprobcycleout}, it is easy to construct an example of a quiver whose complete path algebra $\Lambda$ satisfies the following property: Neither the 2-term silting object $\Lambda$ nor the 2-term silting object $\Lambda[1]$ may be mutated at any of their direct summands.

\subsection*{Acknowledgements} The main ideas of this work were conceived when MK and HUT visited DLF at the University of Padova in May 2025. They are grateful to DLF and the University of Padova for their hospitality and the memorable time. MK is also grateful to Raphael Bennett-Tennenhaus for helpful explanations regarding complete path and string algebras and to Wassilij Gnedin for reading an earlier version of the manuscript and his invaluable suggestions.

\section{Background}

Let $Q$ be a finite quiver, that is, an oriented graph. Let $K$ be a field and denote the path $K$-algebra of $Q$ by $K Q$. Write $\mathfrak{a}$ for the ideal of $KQ$ generated by the arrows of $Q$ and define the completed path algebra $\widehat{KQ}$ as
\[ \widehat{KQ} \coloneqq \lim \limits_{\longleftarrow} KQ/\mfrak^p. \]
Write $KQ_n$ for the $K$-span of paths of length $n$ in $Q$. The elements of $\widehat{KQ}$ are formal infinite sums $\sum_{n \geq 0} a_n$ with $a_n \in KQ_n$ and the product of two elements is given by
\[ \left( \sum_{i \geq 0} a_i \right) \left( \sum_{j \geq 0} b_j \right) = \sum_{k \geq 0} \sum_{i+j = k} a_i b_j. \]
The complete path algebra $\widehat{KQ}$ is a topological $K$-algebra via the $\mathfrak{a}$-adic topology. We fix this setting for the rest of this paper and write $\Lambda = \widehat{KQ}$ throughout. We consider right $\Lambda$-modules and denote by $P_i$ the indecomposable projective $\Lambda$-module (finitely) generated by the idempotent $e_i$ associated to $i \in Q_0$.

The bounded homotopy category of finitely generated projective $\Lambda$-modules $\mathrm{K}^b(\proj \Lambda)$ is a triangulated Krull--Schmidt category \cite[Lem. 2.17]{KellerYang2011} . For $i \in \mathbb{Z}$, we denote the $i$-th shift (also called suspension) of an object $X \in \K^b(\proj \Lambda)$ by $X[i]$.

\begin{definition}(\cite{Aihara_2012})
	An object $X \in \mathrm{K}^b(\proj \Lambda)$ is called \textit{presilting} if $\Hom_{\mathrm{K}^b(\proj \Lambda)}(X,X[i])= 0$ for all $i > 0$ and \textit{silting} if additionally $\thick(X) = \mathrm{K}^b(\proj \Lambda)$, which is to say that $X$ generates $\mathrm{K}^b(\proj \Lambda)$ under taking direct summands, cones and shifts.
\end{definition}

Because $ \mathrm{K}^b(\proj \Lambda)$ is Krull--Schmidt, any silting object has $|Q_0|$ pairwise nonisomorphic indecomposable direct summands \cite[Cor. 2.28]{Aihara_2012}. We call a (pre-)silting object $X \in \K^b(\proj \Lambda)$ \textit{basic} if its indecomposable direct summands are pairwise nonisomorphic and \textit{2-term} if $X$ is isomorphic in $\mathrm{K}^b(\proj \Lambda)$ to a complex concentrated in cohomological degrees $-1$ and $0$. We write $\silt (\Lambda)$ for the set of isomorphism classes of basic silting objects in $\mathrm{K}^b(\proj \Lambda)$ and $\twosilt (\Lambda)$ for the subset of 2-term silting objects therein. 

\begin{theorem}(\cite[Thm. 2.11]{Aihara_2012})
	There exists a partial order $\leq$ on $\silt (\Lambda )$ such that $Y \geq X$ whenever $\Hom_{\mathrm{K}^b(\proj \Lambda)}(Y, X[1]) = 0$. 
\end{theorem}

In some cases, this partial order can be generated via a well-behaved mutation theory.

\begin{theorem}(\cite[Def. 2.30, Thm. 2.31]{Aihara_2012})\label{thm:siltmutdefn}
    Let $X \oplus Y$ be a basic silting object in $\mathrm{K}^b(\proj \Lambda)$. Assume that there exists a left $(\add X)$-approximation $f$ of $Y$ in $\mathrm{K}^b(\proj \Lambda)$. Then $X \oplus \cone(f)$ is a silting complex in $\mathrm{K}^b(\proj \Lambda)$, which we call the \textit{left silting mutation} of $X \oplus Y$ at $Y$.
    
    Dually, assume that there exists a right $(\add X)$-approximation $g$ of $Y$ in $\mathrm{K}^b(\proj \Lambda)$. Then $X \oplus \cone(g)[-1]$ is a silting complex in $\mathrm{K}^b(\proj \Lambda)$, which we call the \textit{right silting mutation} of $X \oplus Y$ at $Y$. 
\end{theorem}

Even if a left or right mutation of a 2-term silting object exists, it may not yield another 2-term silting object (concentrated in degrees $-1$ and $0$). For instance, a left mutation of the silting object $\Lambda[1]$ at any of its indecomposable direct summands, would result in an object supported in degree $-2$. If the summand $Y$ in \cref{thm:siltmutdefn} is indecomposable, then the left/right mutation of $X \oplus Y$ at $Y$ is called \textit{irreducible}. If there does not exist a left $(\add X)$-approximation of $Y$, we say that $X \oplus Y$ is \textit{left immutable} at $Y$ and we define right immutability dually. Moreover, given $X,X' \in \silt (\Lambda)$, write $X \lessdot_{\textnormal{mut}} X'$ if there is an indecomposable direct summand $Y'$ of $X'$ such that $X$ is the left mutation of $X'$ at $Y'$. 

\begin{lemma}\label{lem:submutationposet}
	The reflexive and transitive closure of $\lessdot_{\textnormal{mut}}$ defines a subposet $(\silt (\Lambda), \leq_{\textnormal{mut}})$ of $(\silt (\Lambda), \leq)$. 
\end{lemma} 
\begin{proof}
	Since the underlying sets are the same, we just need to show that $X \lessdot_{\textnormal{mut}} Y$ implies $X \leq Y$, that is, $\Hom_{\mathrm{K}^b(\proj \Lambda)}(Y, X[1]) = 0$. By assumption we can write $X = Y/Y' \oplus \cone(g)$ with $g$ a left $\add (Y/Y')$-approximation of $Y'$. From the assumption that $Y$ is presilting it follows that $\Hom_{\mathrm{K}^b(\proj \Lambda)}(Y, Y/Y'[1]) = 0$. It remains to show that $\Hom_{\mathrm{K}^b(\proj \Lambda)}(Y, \cone(g)[1]) = 0$. This follows by applying the left-exact functor $\Hom_{\mathrm{K}^b(\proj \Lambda)}(Y,-)$ to the triangle 
	\[ Y'[1] \xrightarrow[]{{g[1]}} Y[1] \to \cone(g)[1] \to Y'[2]\]
	 and the fact that $Y$ is presilting. 
\end{proof}

Under certain assumptions, for example when $\add(P)$ is functorially finite for all basic presilting objects $P$ with $|Q_0|-1$ indecomposable direct summands, the subposet $\leq_{\textnormal{mut}}$ coincides with $\leq$, see \cite[Sec. 2.6]{Aihara_2012} and \cite[Sec. 2.4]{Kimura2024}. In contrast to Noetherian $R$-algebras \cite[Lem. 2.3]{Kimura2024}, this assumption is generally not met for complete path algebras (see \cref{eq:sec3}). 

The following example gives a small complete path algebra which is not a Noetherian $R$-algebra for any ring $R$. 

\begin{example}\label{exmp:notNoetherian}
	Let $Q= \begin{tikzcd} \bullet \arrow[loop, in=60, out=-60, looseness=5, "x",swap] \arrow[loop, in=120, out=240, looseness=5, "y"]\end{tikzcd}$, then $KQ = K\langle x,y \rangle$ and $\widehat{KQ} = K \llangle x,y \rrangle$. The $K$-algebra $\widehat{KQ}$ is not finitely generated over $K$ and is neither itself commutative nor Noetherian. Hence it cannot be a Noetherian $R$-algebra in the sense of \cite{Ausl1977}.
\end{example}

\section{Proofs of the main results}\label{eq:sec3}

In this section, we establish certain connected subquivers $Q'$ of $Q$ which pose obstructions to the existence of the desired isomorphisms of partially ordered sets. Indeed, from these subquivers we construct silting objects which are immutable at (at least one) on of their indecomposable direct summands. These form one of the key ingredients in the proof \cref{thm:main1}. 

For two distinct elements $x$ and $y$ of a poset $(P,\leq)$, write $x \lessdot y$ if $x \leq y$ and $x \leq z \leq y$ implies $z \in \{x,y\}$. The proof of \cref{thm:main1} uses the Hasse quiver of a poset $P$, which is defined as the quiver whose elements are elements of $P$ and for which there exists an arrow $y \leftarrow x$ whenever $y \lessdot x$. It is well-known and easy to show that a poset isomorphism $f: P \to P'$ gives rise to an isomorphism of Hasse quivers. In particular, we have $x \lessdot y$ if and only if $f(x) \lessdot f(y)$. 

In the following lemmas we view projective $\Lambda$-modules as stalk complexes sitting in degree 0 in $\mathrm{K}^b(\proj \Lambda)$.

\begin{lemma}\label{lem:mutprobcyclein}
	Let $Q$ be a quiver and let $\Lambda = \widehat{KQ}$. If $Q$ contains a subquiver of the form 
	\[ \begin{tikzcd} 
	& v_{r-1} \arrow[rd, bend left, "\epsilon_r"] \\
	\vdots \arrow[ru, bend left, "\epsilon_{r-1}"] &&v \arrow[ld, bend left, "\epsilon_1"] & w \arrow[l, "\alpha"] \\
	&v_1 \arrow[lu, bend left, "\epsilon_2"]
	 \end{tikzcd}, \]
	with $w \not \in \{v, v_1, \dots, v_{r-1}\}$, then the 2-term silting object $\Lambda$ in $\mathrm{K}^b(\proj \Lambda)$ is left immutable at $P_v$. 
\end{lemma}
\begin{proof}
	Let $\epsilon = \epsilon_1 \dots \epsilon_r$, and note that possibly $r=1$, in which case $\epsilon$ is a loop. By \cite[Lem. 2.2]{Kimura2024}, it is equivalent to show that $\Hom_{\Lambda}(P_v, \Lambda/P_v)$ is an infinitely generated left $\End_{\Lambda}(\Lambda/P_v)$-module. Whence it is sufficient to show that $\Hom_{\Lambda}(P_v, P_w)$ is an infinitely generated left $\End_{\Lambda}(\Lambda/P_v)$-module. Because $\Hom_{\Lambda}(P(v), P(w)) = e_w \Lambda e_v$, it is sufficient to show that there are infinitely many paths $\{\gamma_i\}_{i \in \mathcal{I}}$ from $w$ to $v$, which are not of the form $\gamma_i = c_i \gamma_i'$ with $c_i$ a path from $w$ to a vertex different to $v$. It is immediate that the infinite set $\{\alpha \epsilon^r\}_{r\geq0}$ is such a set.
\end{proof}

We also record the dual result.

\begin{lemma}\label{lem:mutprobcycleout}
	Let $Q$ be a quiver and let $\Lambda = \widehat{KQ}$. If $Q$ contains a subquiver of the form 
	\[ \begin{tikzcd} 
	& v_{r-1} \arrow[rd, bend left, "\epsilon_r"] \\
	\vdots \arrow[ru, bend left, "\epsilon_{r-1}"] &&v \arrow[ld, bend left, "\epsilon_1"] \arrow[r, "\alpha"] & w \\
	&v_1 \arrow[lu, bend left, "\epsilon_2"]
	 \end{tikzcd}, \]
	with $w \not \in \{ v, v_1, \dots, v_{r-1}\}$, then the 2-term silting object $\Lambda[1]$ is right immutable at $P_v[1]$. 
\end{lemma}

Finally, we need to consider the following case, where an indecomposable direct summand only become immutable after an initial mutation.

\begin{lemma}\label{lem:mutprobAtilde}
	Let $Q$ be a quiver containing a subquiver of the form 
	\[ \begin{tikzcd}  v \arrow[r, "\alpha"] \arrow[r, bend left=50, "\epsilon"] & w \arrow[l, bend left=50, "\beta"] \end{tikzcd} \]
	with $w\neq v$, and let $\Lambda = \widehat{KQ}$. Assume that a left mutation of $\Lambda$ at $P_w$ exists and denote it by $\mu_{P_w}(\Lambda)$. Then the 2-term silting complex $\mu_{P_w}(\Lambda)$ is left immutable at $P_v$. 
\end{lemma}
\begin{proof}
	Assume that $\mu_{P_w}(\Lambda)$ exists, which is to say that there exists a left $\add(\Lambda/P_w)$-approximation $g$ of $P_w$. Such an approximation has component morphisms determined by incoming arrows at $w$, in particular by $\alpha$ and $\epsilon$. For this reason, we may write $g$ as
	\[
	P_w \xrightarrow{(\alpha \phantom{x} \epsilon \phantom{x} g')^t} P_v \oplus P_v \oplus Q,
	\]
	where $Q$ is a direct sum of indecomposable projective $\Lambda$-modules different from $P_w$. Denote by $X = \cone(g) \in K^{[-1,0]}(\proj \Lambda)$ the corresponding 2-term presilting object. The 2-term silting object $\mu_{P_w}(\Lambda)$ is then given by $X \oplus (\Lambda/P_w)$. We now show that $\mu_{P_w}(\Lambda)$ is left immutable at $P_v$. Assume for a contradiction that there exists a left $\add(\mu_{P_w}(\Lambda)/P_v)$-approximation $h$ of $P_v$ and write $Y$ for the codomain of $h$. 
	
	Observe that for each pair $(i,j) \in \mathbb{Z}_{\geq 1}^2$, there exists a morphism $\phi_{ij}: P_v[0] \to X$ of 2-term complexes given by
	\[ 
	\begin{tikzcd}[ampersand replacement=\&, column sep=50, row sep = 50]
	0 \arrow[r, ] \arrow[d] \& P_v \arrow[d, "\begin{pmatrix} (\epsilon \beta)^i \\ (\alpha \beta)^j \\ 0\end{pmatrix}"] \\
	P_w \arrow[r, "{\begin{pmatrix} \alpha & \epsilon & g' \end{pmatrix}}^t"] \& P_v \oplus P_v \oplus Q
	\end{tikzcd}.
	\] 
	Since $h$ is a left $\add(\mu_{P_w}(\Lambda)/P_v)$-approximation, there exists a map $p_{ij}: Y \to X$ for each $(i,j) \in \mathbb{Z}_{\geq 1}^2$ such that $p_{ij}h = \phi_{ij}$. There are two ways of writing the component morphism $(\epsilon \beta)^i$ of $\phi_{ij}$ as a composite of two morphisms. Decomposing $(\epsilon \beta)^i=  ((\epsilon \beta)^{i-s-1})\epsilon) \circ (\beta(\epsilon \beta)^{s})$, for some $s \geq 1$, yields a contradiction to $Y \in \add(\mu_{P_w}(\Lambda)/P_v)$. Indeed, the projective $\Lambda$-module which forms the degree zero part of $Y$ would need to contain $P_w$ as a direct summand because $P_w$ is the codomain of $\beta(\epsilon \beta)^s$. Decomposing $(\epsilon \beta)^i = (\epsilon \beta)^{i-s} \circ (\epsilon \beta)^s$, for some $s \geq 1$, implies that the projective $\Lambda$-module which forms the degree zero part of $Y$ contains $P_v$ as a direct summand because $P_v$ is the codomain of $(\epsilon \beta)^s$. Since $Y \in \add(\mu_{P_w}(\Lambda)/P_v)$, this implies that $X$ is a direct summand of $Y$ and that the chain map $h: P_v[0] \to Y$ has a component chain map $h': P_v[0] \to X$. Moreover, the degree zero map of $h'$ has a component morphism given by $(\epsilon \beta)^s: P_v \to P_v$, for some $s \geq 1$. Similarly, we deduce from the decompositions of $(\alpha \beta)^j$, that $h$ has a component chain map $h'': P_v[0] \to X$ whose degree zero map has a component morphism given by $(\alpha \beta)^r: P_v \to P_v$, for some $r \geq 1$. Therefore, without loss of generality, we may assume that $h$ has a component chain map $h'$ of the one of following four forms:
	\[ 
	\begin{tikzcd}[ampersand replacement=\&, column sep=50, row sep = 50]
	0 \arrow[r, ] \arrow[d] \& P_v \arrow[d, "\begin{pmatrix} (\epsilon \beta)^s \\ (\alpha \beta)^r \\ q \end{pmatrix}"] \\
	P_w \arrow[r, "{\begin{pmatrix} \alpha & \epsilon & g' \end{pmatrix}}^t"] \& P_v \oplus P_v \oplus Q 
	\end{tikzcd}
	\qquad  \qquad
	\begin{tikzcd}[ampersand replacement=\&, column sep=50, row sep = 50]
	0 \arrow[r, ] \arrow[d] \& P_v \arrow[d, "\begin{pmatrix} (\alpha \beta)^r \\ (\epsilon \beta)^s \\ q \end{pmatrix}"] \\
	P_w \arrow[r, "{\begin{pmatrix} \alpha & \epsilon & g' \end{pmatrix}}^t"] \& P_v \oplus P_v \oplus Q 
	\end{tikzcd}.
	\] 	
	\[ 
	\begin{tikzcd}[ampersand replacement=\&, column sep=50, row sep = 80]
	0 \arrow[r, ] \arrow[d] \& P_v \arrow[d, "\begin{pmatrix} (\epsilon \beta)^s \\ q_1 \\ q_2 \\ q_3 \\ (\alpha \beta)^r \\ q_4  \end{pmatrix}"] \\
	\big( P_w \arrow[r, "{\begin{pmatrix} \alpha & \epsilon & g' \end{pmatrix}}^t"] \& P_v \oplus P_v \oplus Q \big)^2
	\end{tikzcd}
	\qquad  \qquad
	\begin{tikzcd}[ampersand replacement=\&, column sep=50, row sep = 80]
	0 \arrow[r, ] \arrow[d] \& P_v \arrow[d, "\begin{pmatrix} (\alpha \beta)^r \\ q_1 \\ q_2 \\ q_3 \\ (\epsilon \beta)^s \\ q_4  \end{pmatrix}"] \\
	\big( P_w \arrow[r, "{\begin{pmatrix} \alpha & \epsilon & g' \end{pmatrix}}^t"] \& P_v \oplus P_v \oplus Q \big)^2
	\end{tikzcd}.
	\] 
	Since $\phi_{ij}: P_v[0] \to X$ factors through $h'$, there is an endomorphism $x \in \End(X)$, which forms a component chain map of $p_{ij}$ such that $x h' = \phi_{ij}$. We show that such an $x$ can only exist if $s=i$ and $r=j$. We only address the top two cases in order from left to right. The bottom two cases are analogous. 
	
	In the top-left case, the degree zero map of $x$ needs to be of the form 
	\[ x_0 = \begin{pmatrix} (\epsilon \beta)^{i-s} & * & * \\ * & (\alpha \beta)^{j-r} & * \\ * & * &*\end{pmatrix}. \]
	If $s<i$ or $r <j$, a direct calculations shows that there does not exist an endomorphism $x_{-1} \in \End(P_w)$ making the necessary square (to make $x$ into a chain map) commute. Similarly, in the top-right case, the degree zero map of $x$ needs to be of the form 
	\[ x_0 = \begin{pmatrix} * & (\epsilon \beta)^{i-s} & * \\ (\alpha \beta)^{j-r} &* & * \\ * & * &*\end{pmatrix}. \]
	In this case, one sees that there cannot exist an endomorphism $x_{-1} \in \End(P_w)$ making the necessary square commute. The bottom two cases follow similarly. 
	
	This shows that each $\phi_{ij}: P_v[0] \to X$ has to be a component of the left $\add(\mu_{P_w}(\Lambda)/P_v)$-approximation $h$. Since there are infinitely many such maps, the codomain of $h$ would have to be infinite dimensional, yielding a contradiction. This concludes the proof. 
	\end{proof}

In the following proof of \cref{thm:main1}, we write $|\Lambda|$ for the number of nonisomorphic simple $\Lambda$-modules, which equals $|Q_0|$. We call $\Lambda$ local if it has a unique maximal left (equiv. right) ideal. 

\begin{proof}[Proof of \cref{thm:main1}]
    $(1) \Rightarrow (3)$ The isomorphism $f$ induces an isomorphism of Hasse quivers of the two posets. Combining \cite[Cor. 4.6]{dirrt} and \cite[Cor. 2.34, Cor. 3.9]{tau}, it follows that the graph underlying the Hasse quiver $\Hasse(\twosilt, \leq)$ is a $|\Lambda/I|$-regular graph, meaning that there are $|\Lambda/I|$ edges incident to every vertex. Therefore, the graph underlying the Hasse quiver of $(\twosilt \Lambda, \leq_{\textnormal{mut}})$ is also $|\Lambda/I|$-regular. Since $I \subseteq \rad \Lambda$, no primitive idempotent $e \in \Lambda$ is contained in $I$, whence $|\Lambda/I| = |\Lambda|$ holds. Consequently, the Hasse quiver of $(\twosilt \Lambda, \leq_{\textnormal{mut}})$ is $|\Lambda|$-regular. 
    
    As a consequence of \cite[Lem. 2.12]{Kimura2024}, no 2-term silting object in $\mathrm{K}^b(\proj \Lambda)$ is both left and right mutable at any of its indecomposable direct summands. However, since every basic 2-term silting object over $\Lambda$ has $|\Lambda|$ indecomposable direct summands \cite[Thm. 2.27, Cor. 2.28]{Aihara_2012}, this implies that every basic 2-term silting object over $\Lambda$ is either left or right mutable at each of its indecomposable direct summands.
    
    Assume that $Q$ is not acyclic, has more than one vertex and is not just an oriented cycle. In this case, $Q$ contains a connected subquiver $Q'$ as in \cref{lem:mutprobcyclein}, \cref{lem:mutprobcycleout} or \cref{lem:mutprobAtilde}. If $Q'$ is as in \cref{lem:mutprobcyclein}, the lemma states that the 2-term silting stalk complex $\Lambda$ concentrated in degree zero is left immutable at one of its indecomposable direct summands. If a right mutation at any indecomposable direct summand of $\Lambda$ existed, the resulting object would be concentrated in degrees $0$ and $1$. In other words, it would not be a 2-term silting object. It follows that there exists a 2-term silting object with a left and right immutable indecomposable direct summand, which contradicts the conclusion of the previous paragraph. If $Q'$ is as in \cref{lem:mutprobcycleout}, the statement of the lemma, and a dual argument about degrees as above, leads to the same contradiction. 
    
     If $Q'$ is as in \cref{lem:mutprobAtilde}, then there are two cases: If there already exists an indecomposable direct summand of the stalk complex $\Lambda$ concentrated in degree zero which is left immutable, then the argument of the previous paragraph yields a contradiction. Otherwise, \cref{lem:mutprobAtilde} establishes the existence of a 2-term silting object $X$ which is left immutable at one of its indecomposable direct summands. Since this indecomposable direct summand is concentrated in degree zero, the argument of the previous paragraph shows that $X$ is also left immutable at this indecomposable direct summand, yielding a contradiction.
     
     $(2) \Rightarrow (3)$ This implication is similar to the one above. The isomorphism $f'$ induces an isomorphism of Hasse quivers of the two posets. For the finite dimensional quotient algebra $\Lambda/I$, the poset $(\silt \Lambda/I, \leq)$ is $(2 | \Lambda| )$-regular because the required approximations always exist as $\Hom$-spaces are finite dimensional. This means that the graph underlying $\Hasse(\silt \Lambda, \leq_{\textnormal{mut}})$ is $(2 |\Lambda|)$-regular, or in other words, that each silting object is both left and right mutable at each of its indecomposable direct summands. If $Q$ is not one of the three classes of quivers we want to obtain, it contains a subquiver as in \cref{lem:mutprobcyclein}, \cref{lem:mutprobcycleout} or \cref{lem:mutprobAtilde}. In each case, the lemma yields the existence of a silting object, which cannot be mutated at one of its indecomposable direct summands, yielding a contradiction.
     
     $(3) \Rightarrow (4)$ If $Q$ has only one vertex, it has a unique maximal ideal generated by its arrows, making $\Lambda$ a local algebra. If $Q$ is acyclic, then $\Lambda$ is a finite dimensional $K$-algebra. In particular, it is a finitely generated $K$-module. If $Q$ is an oriented cycle, let $x$ be the sum of all cyclic paths of minimal length in the quiver and let $R = K\llbracket x \rrbracket$. Then a standard argument shows that $\Lambda$ is an $R$-algebra and finitely generated as an $R$-module, see \cite[Thm. 3.2.8]{Ricke}. This is a well-known hereditary $K \llbracket x \rrbracket$-order. 
     
     $(4) \Rightarrow (1)$ If $\Lambda$ is a local algebra (finite dimensional or not), then it has only one nonzero idempotent $e \in \Lambda$ by \cite[Prop. 19.2]{Lam2001}. Thus, up to isomorphism, there exists exactly one indecomposable projective $\Lambda$-module $P$. This means that $\K^b(\proj \Lambda)$ has an indecomposable silting object. It follows that the Hasse quiver of 2-term silting objects is of the form $P[0] \to P[1]$, see \cite[Thm. 2.26]{Aihara_2012}. In the other two cases, $\Lambda$ is a Noetherian $R$-algebra with $R$ a commutative complete local Noetherian ring $R$. Let $\mathfrak{m} \subseteq R$ be the maximal ideal of $R$. By \cite[Thm. 4.3]{Kimura2024}, any ideal $I \subseteq \mathfrak{m} \Lambda$ has the desired property. By \cite[Prop. 5.22]{CurtisReiner1981}, we have $\mathfrak{m} \Lambda \subseteq \rad \Lambda$, which completes the proof. 
     
     $(4) \Rightarrow (2)$ This implication is similar to the previous one. Using \cite[Thm. 2.26]{Aihara_2012}, we deduce that if $\Lambda$ is a local algebra, then $\Hasse(\silt \Lambda, \leq)$ is given by
     \[ \dots \Lambda[-1] \to \Lambda[0] \to \Lambda[1] \to \dots \]
     when viewing $\Lambda$ as an indecomposable stalk complex in $\mathrm{K}^b(\proj \Lambda)$. The implication follows. If $\Lambda$ is acyclic, the implication holds trivially for the ideal $I= 0$. If $\Lambda$ is a hereditary $K \llbracket x \rrbracket$-order, let $\mathfrak{n}$ be the maximal ideal of $R$. Then the implication holds for the ideal generated by $\mathfrak{n}$, by \cite[Thm. A]{Gnedin2022}, because the ring $K \llbracket x \rrbracket$ is complete regular local. Again, by \cite[Prop. 5.22]{CurtisReiner1981}, we have $\mathfrak{n} \Lambda \subseteq \rad \Lambda$, as required.
        \end{proof}
        
        \cref{cor:main1} is a direct consequence of this.

    \begin{proof}[Proof of \cref{cor:main1}]
    	If $\Lambda$ is a module-finite $R$-algebras with $R$ as in \eqref{cor:main1eq1}, then \cite[Thm. 4.3]{Kimura2024} implies \cref{thm:main1}\eqref{thm:main1eq1}. By \cref{thm:main1}\eqref{thm:main1eq1}$\Rightarrow$\eqref{thm:main1eq3}, this means that $Q$ has one vertex or meets the conditions in \eqref{cor:main1eq2}. By assumption, $Q$ has more than one vertex. The converse is \cref{thm:main1}\eqref{thm:main1eq3}$\Rightarrow$\eqref{thm:main1eq4}.
    \end{proof}

We continue with a proof of the characterisation of quivers for which all finite dimensional quotients are 2-term silting finite. In the following, we freely use that 2-term silting finiteness of finite dimensional $K$-algebras is equivalent to $\tau$-tilting finiteness \cite[Thm. 0.5]{tau} without introducing the latter concept. Roughly speaking, this translation replaces 2-term silting complexes by their $0$-th homologies, which are right $\Lambda$-modules. The most important and obvious consequence of this substitution is that representation finite algebras are $\tau$-tilting finite and hence 2-term silting finite.

\begin{proof}[Proof of \cref{thm:main2}]
$(1) \Rightarrow (2)$ is immediate. 

$(2) \Rightarrow (3)$ Assume that $Q$ is not of Dynkin type ADE, has more than one vertex and is not an oriented cycle. 
If $Q$ is acyclic then for $i > |Q_1|$ the algebra $\Lambda/\mathfrak{a}^i \simeq KQ$ is hereditary. Since $Q$ is not of finite Dynkin type, $\Lambda/\mathfrak{a}^i $ is representation infinite. A hereditary algebra is representation infinite if and only if it is $\tau$-tilting infinite as a consequence of \cite[Cor. IV.4.7]{auslander_reiten_smalo_1995}. This is a contradiction.

If $Q$ contains a cycle $\epsilon_1 \dots \epsilon_r$, but is not just an oriented cycle, it contains a subquiver of the form as in \cref{lem:mutprobcyclein} or \cref{lem:mutprobcycleout} with $w \neq v$ because $Q$ has more than one vertex, but possibly $w \in \{ v_1, \dots, v_{r-1}\}$ if $r \geq 2$. The orientation of $\alpha$ may be arbitrary in the following. Note that for $r=1$, the arrow $\epsilon_1$ is a loop. We distinguish between two further cases.

Assume $r \geq 2$ and $w \in \{ v_1, \dots, v_{r-1}\}$. For $i > |Q_1|$, the algebra $\Lambda/\mathfrak{a}^i$ admits a quotient algebra of hereditary type $\widetilde{A}$ in this case. Since this hereditary quotient is representation infinite, it is also $\tau$-tilting infinite as before. From \cite[Cor. 1.9]{dirrt} we deduce that $\Lambda/\mathfrak{a}^i$ is not $\tau$-tilting finite in this case, which yields a contradiction.

Assume $w \not \in \{ v_1, \dots, v_{r-1}\}$. Let $e_{v}$ and $e_w$ be the primitive idempotents corresponding to the vertices $v$ and $w$ of $Q$ respectively and let $i \geq 4|Q_1|$. Then the algebra $\Gamma \coloneqq (e_{v} + e_w) (\Lambda/\mathfrak{a}^i \Lambda) (e_{v} + e_w)$ has a quotient algebra $\Gamma/I$ isomorphic to 
\[
K \Bigg( \begin{tikzcd}
    v \arrow[loop, "\epsilon", in=135, out = 225, looseness=5] \arrow[r, dash, "\alpha"] & w 
\end{tikzcd} \Bigg) / \langle \epsilon^4 \rangle .
\]

The algebra $\Gamma/J$ is $\tau$-tilting infinite regardless of the orientation of $\alpha$ by \cite[Thm. 1.1]{Wang2025} and \cite[Thm. 2.14]{tau}. Thus, $\Gamma$ is $\tau$-tilting infinite by \cite[Cor. 1.9]{dirrt}. From \cite[Thm. 1.4]{dij17}, it follows that there exist infinitely many nonisomorphic bricks in $\mods \Gamma$. By \cite[Thm. I.6.8]{assem_skowronski_simson_2006}, there is a fully faithful functor $F: \mods \Gamma \to \mods (\Lambda/\mathfrak{a}^i)$ which preserves bricks \cite[Lem. 2.1]{Zhang2023}. Therefore, there are finitely many nonisomorphic bricks in $\mods (\Lambda/\mathfrak{a}^i)$, which implies that $\Lambda/\mathfrak{a}^i$ is $\tau$-tilting infinite \cite[Thm. 1.4]{dij17}. This is a contradiction.

Note that $r=1$ and $w \in \{ v_1, \dots, v_{r-1}\}$ is impossible because $Q$ has more than one vertex. This completes the proof of this implication.

$(3) \Rightarrow (1)$ If $Q$ is of Dynkin type ADE, then $\Lambda \simeq KQ$ holds because the underlying quiver is acyclic. Thus, $\Lambda$ is finite dimensional representation finite as it is of finite Dynkin type. Consequently, all quotient algebras are representation finite and thus $\tau$-tilting finite. Similarly, if $Q$ is an oriented cycle, then any finite dimensional quotient $\Lambda/I$ is a so-called Nakayama algebra, which are representation finite \cite[Thm. V.3.5]{assem_skowronski_simson_2006} and thus $\tau$-tilting finite. If $Q$ has only one vertex, then for any finite dimensional quotient $\Lambda/I$, the zero module is an almost complete $\tau$-tilting pair. By \cite[Thm. 0.4]{tau} such a pair has exactly two completions. This gives a finite connected component of the support $\tau$-tilting quiver of $\Lambda/I$, which means that the algebra $\Lambda/I$ is $\tau$-tilting finite \cite[Cor. 2.38]{tau}. 
\end{proof}

\subsection*{AI contribution disclosure}
The authors did not use Artificial Intelligence to any extent in the research process or during manuscript preparation.

\printbibliography

@article{tau, 
title={$\tau$-tilting theory}, 
volume={150}, 
number={3}, 
journal={Compos. Math.}, 
publisher={London Mathematical Society}, author={Adachi, Takahide and Iyama, Osamu and Reiten, Idun}, 
year={2014},
pages={415–452}}

@article {DerksenFei2015,
    AUTHOR = {Derksen, Harm and Fei, Jiarui},
     TITLE = {General presentations of algebras},
   JOURNAL = {Adv. Math.},
  FJOURNAL = {Advances in Mathematics},
    VOLUME = {278},
      YEAR = {2015},
     PAGES = {210--237},
}

@article{Wang2025,
 author = {Wang, Qi},
 title = {$\tau$-tilting finiteness of two-point algebras {II}},
 fjournal = {Journal of Algebra and its Applications},
 journal = {J. Algebra Appl.},
 volume = {24},
 number = {2},
 pages = {33},
 year = {2025},
}

@book{assem_skowronski_simson_2006, place={Cambridge}, series={London Mathematical Society Student Texts}, title={Elements of the Representation Theory of Associative Algebras: Techniques of Representation Theory}, volume={1}, publisher={Cambridge University Press}, author={Assem, Ibrahim and Skowronski, Andrzej and Simson, Daniel}, year={2006}, collection={Lond. Math. Soc. Stud. Texts}}

@article{dij17,
	author = {Demonet, Laurent and Iyama, Osamu and Jasso, Gustavo},
	title = {$\tau$-tilting finite algebras, bricks, and $g$-vectors},
	journal = {Int. Math. Res. Not.},
	volume = {2019},
	number = {3},
	pages = {852-892},
	year = {2017},
}

@book{auslander_reiten_smalo_1995, place={Cambridge}, series={Camb. Stud. Adv. Math.}, title={Representation Theory of Artin Algebras}, publisher={Cambridge University Press}, author={Auslander, Maurice and Reiten, Idun and Smal{\o}, Sverre O.}, year={1995}, collection={Camb. Stud. Adv. Math.}}

@article{dirrt,
  author = {Demonet, Laurent and Iyama, Osamu and Reading, Nathan and Reiten, Idun and Thomas, Hugh},
  title = {Lattice theory of torsion classes: Beyond $\tau$-tilting theory},
    year = {2023},
    month = {04},
    pages = {542-612},
    volume = {10},
    journal = {Trans. Am. Math. Soc., Ser. B},
}

@article {BMRRT,
    AUTHOR = {Buan, Aslak B. and Marsh, Bethany R. and Reineke, Markus and
              Reiten, Idun and Todorov, Gordana},
     TITLE = {Tilting theory and cluster combinatorics},
   JOURNAL = {Adv. Math.},
  FJOURNAL = {Advances in Mathematics},
    VOLUME = {204},
      YEAR = {2006},
    NUMBER = {2},
     PAGES = {572--618},
}

@article{Aihara_2012,
   title={Silting mutation in triangulated categories},
   volume={85},
   number={3},
   journal={J. Lond. Math. Soc., II. Ser.},
   author={Aihara, Takuma and Iyama, Osamu},
   year={2012},
   }

@article{Zhang2023,
title = {A construction of support $\tau$-tilting modules over $\tau$-tilting finite algebras},
author = {Zhang, Yingying},
journal = {Adv. Math., Beijing},
volume = {53},
pages = {1039-1047},
year = {2023},
}

@article{KellerYang2011,
 author = {Keller, Bernhard and Yang, Dong},
 title = {Derived equivalences from mutations of quivers with potential},
 fjournal = {Advances in Mathematics},
 journal = {Adv. Math.},
 volume = {226},
 number = {3},
 pages = {2118--2168},
 year = {2011},
}

@article{Kimura2024,
 author = {Kimura, Yuta},
 title = {Tilting and silting theory of {Noetherian} algebras},
 fjournal = {IMRN. International Mathematics Research Notices},
 journal = {Int. Math. Res. Not.},
 issn = {1073-7928},
 volume = {2024},
 number = {2},
 pages = {1685--1732},
 year = {2024},
}

@article{DWZ2008,
 author = {Derksen, Harm and Weyman, Jerzy and Zelevinsky, Andrei},
 title = {Quivers with potentials and their representations. {I}: {Mutations}.},
 fjournal = {Selecta Mathematica. New Series},
 journal = {Sel. Math., New Ser.},
 volume = {14},
 number = {1},
 pages = {59--119},
 year = {2008},
}

@article{DWZ2010,
 author = {Derksen, Harm and Weyman, Jerzy and Zelevinsky, Andrei},
 title = {Quivers with potentials and their representations. {II}: {Applications} to cluster algebras.},
 fjournal = {Journal of the American Mathematical Society},
 journal = {J. Am. Math. Soc.},
 volume = {23},
 number = {3},
 pages = {749--790},
 year = {2010},
}

@phdthesis{Ricke,
author = {Charlotte Ricke},
title = {On $\tau$-tilting theory and perpendicular categories},
school = {Rheinische Friedrich-Wilhelms-Universität Bonn},
year = 2016,
month = nov,
}

@misc{Ausl1977,
 author = {Auslander, Maurice},
 title = {Existence theorems for almost split sequences},
 year = {1977},
 howpublished = {Ring {Theory} {II}, {Proc}. 2nd {Okla}. {Conf}. 1975, 1-44.},
}

@misc{Gnedin2022,
      title={Silting theory under change of rings}, 
      author={Wassilij Gnedin},
      year={2022},
      eprint={2204.00608},
      archivePrefix={arXiv},
      primaryClass={math.RT},
}

@article{BT2025,
 author = {Bennett-Tennenhaus, Raphael},
 title = {String algebras over local rings: admissibility and biseriality},
 fjournal = {Journal of Algebra},
 journal = {J. Algebra},
 volume = {667},
 pages = {325--364},
 year = {2025},
}

@misc{BT2024,
      title={String algebras over local rings: regular examples}, 
      author={Raphael Bennett-Tennenhaus},
      year={2024},
      eprint={2405.02702},
      archivePrefix={arXiv},
      primaryClass={math.RT},
}

@article{IyamaKimura2024,
 author = {Iyama, Osamu and Kimura, Yuta},
 title = {Classifying subcategories of modules over {Noetherian} algebras},
 fjournal = {Advances in Mathematics},
 journal = {Adv. Math.},
 volume = {446},
 note = {Id/No 109631},
 year = {2024},
}

@misc{CurtisReiner1981,
 author = {Curtis, Charles W. and Reiner, Irving},
 title = {Methods of representation theory, with applications to finite groups and orders. {Vol}. {I}},
 year = {1981},
 howpublished = {Pure and {Applied} {Mathematics}. {A} {Wiley}-{Interscience} {Publication}. {New} {York} etc.: {John} {Wiley} \& {Sons}. {XXI}},
}

@book {Lam2001,
    AUTHOR = {Lam, Tsit-Yuen},
     TITLE = {A first course in noncommutative rings},
    SERIES = {Grad. Texts Math.},
    VOLUME = {131},
   EDITION = {Second},
 PUBLISHER = {Springer},
      YEAR = {2001},
     PAGES = {xx+385},
}

@incollection {KV88,
	AUTHOR = {Keller, Bernhard and Vossieck, Dieter},
	TITLE = {Aisles in derived categories},
	NOTE = {Deuxi\`eme Contact Franco-Belge en Alg\`ebre (Faulx-les-Tombes,
		1987)},
	JOURNAL = {Bull. Soc. Math. Belg. S\'{e}r. A},
	FJOURNAL = {Bulletin de la Soci\'{e}t\'{e} Math\'{e}matique de Belgique. S\'{e}rie A},
	VOLUME = {40},
	YEAR = {1988},
	NUMBER = {2},
	PAGES = {239--253},
}

@article {KY14,
	AUTHOR = {Koenig, Steffen and Yang, Dong},
	TITLE = {Silting objects, simple-minded collections, {$t$}-structures
		and co-{$t$}-structures for finite-dimensional algebras},
	JOURNAL = {Doc. Math.},
	FJOURNAL = {Documenta Mathematica},
	VOLUME = {19},
	YEAR = {2014},
	PAGES = {403--438},
}

@article{Brumer1966,
 author = {Brumer, Armand},
 title = {Pseudocompact algebras, profinite groups and class formations},
 fjournal = {Journal of Algebra},
 journal = {J. Algebra},
 volume = {4},
 pages = {442--470},
 year = {1966},
}

@article{EJR2018,
 author = {Eisele, Florian and Janssens, Geoffrey and Raedschelders, Theo},
 title = {A reduction theorem for $\tau$-rigid modules},
 fjournal = {Mathematische Zeitschrift},
 journal = {Math. Z.},
 volume = {290},
 number = {3-4},
 pages = {1377--1413},
 year = {2018},
}

@article{Gabriel1972,
 author = {Gabriel, Peter},
 title = {Unzerlegbare {Darstellungen}. {I}. ({Indecomposable} representations. {I})},
 fjournal = {Manuscripta Mathematica},
 journal = {Manuscr. Math.},
 volume = {6},
 pages = {71--103},
 year = {1972},
}

@article{Rickard1989,
 author = {Rickard, Jeremy},
 title = {Morita theory for derived categories},
 fjournal = {Journal of the London Mathematical Society. Second Series},
 journal = {J. Lond. Math. Soc., II. Ser.},
 volume = {39},
 number = {3},
 pages = {436--456},
 year = {1989},
}

\end{document}